\documentclass[11pt]{amsart}%
\usepackage{amsmath}
\usepackage{graphicx}%
\usepackage{amsfonts}%
\usepackage{amssymb}
\newtheorem{theorem}{Theorem}
\theoremstyle{plain}

\newtheorem{corollary}[theorem]{Corollary}

\newtheorem{lemma}[theorem]{Lemma}

\newtheorem{proposition}[theorem]{Proposition}

\newcommand{\N}{{\mathbb{N}}}
\newcommand{\G}{{\mathcal{G}}}
\newcommand{\F}{{\mathcal{F}}}
\newcommand{\cH}{{\mathcal{H}}}
\newcommand{\cS}{{\mathcal{S}}}
\newcommand{\A}{{\mathcal{A}}}

\newcommand{\T}{{\mathcal{T}}}
\newcommand{\E}{{\mathcal{E}}}
\newcommand{\D}{{\mathcal{D}}}
\newcommand{\io}{\iota}

\newcommand{\co}{\operatorname{co}}
\newcommand{\supp}{\operatorname{supp}}
\newcommand{\ep}{\varepsilon}
\newcommand{\MF}{\!{}^M\!\F}

\begin{document}
\title{$\ell^{1}$-spreading models in subspaces of mixed Tsirelson spaces}
\author{Denny H. Leung}
\address{Department of Mathematics, National University of Singapore, 2 Science Drive
2, Singapore 117543.}
\email{matlhh@nus.edu.sg}
\author{Wee-Kee Tang}
\address{Mathematics and Mathematics Education, National Institute of Education \\
Nanyang Technological University, 1 Nanyang Walk, Singapore 637616.}
\email{wktang@nie.edu.sg}

\begin{abstract}
We investigate the existence of higher order $\ell^1$-spreading models in subspaces of mixed Tsirelson spaces.  For instance, we show that the following conditions are equivalent for the mixed Tsirelson space $X = T[(\theta_n,\cS_n)^\infty_{n=1}]$: 
\begin{enumerate}
\item Every block subspace of $X$ contains an $\ell^1$-$\cS_{\omega}$-spreading model,
\item The Bourgain $\ell^1$-index $I_b(Y) = I(Y) > \omega^{\omega}$ for any block subspace $Y$ of $X$,
\item $\lim_m\limsup_n\theta_{m+n}/\theta_n > 0$ and every block subspace $Y$ of $X$ contains a block sequence equivalent to a subsequence of the unit vector basis of $X$.
\end{enumerate}
Moreover, if one (and hence all) of these conditions holds, then $X$ is arbitrarily distortable.
\end{abstract}
\maketitle

\section{Introduction}

The discovery and construction of non-trivial asymptotic $\ell^1$ spaces has led to much progress in the structure theory of Banach spaces.  The first such space discovered was Tsirelson's space \cite{T}.  Subsequently, Schlumprecht constructed what is now called Schlumprecht's space \cite{S}.  This space plays a vital role in the solutions of the unconditional basic sequence problem by Gowers and Maurey \cite{GM} and the distortion problem by Odell and Schlumprecht \cite{OS}.  Argyros and Deliyanni \cite{AD} introduced the class of mixed Tsirelson spaces which provides a general framework for Tsirelson's space, Schlumprecht's space and related examples such as Tzafriri's space \cite{Tz}.  Mixed Tsirelson spaces have been studied extensively.  In particular, results about their finite dimensional $\ell^1$-structure were obtained in \cite{ADKM, ADM, M}.
The present authors computed the Bourgain $\ell^1$-indices of mixed Tsirelson spaces in \cite{LT}, and investigated thoroughly the existence of higher order $\ell^1$-spreading models in such spaces \cite{LT2}. (Results in this direction for certain mixed Tsirelson spaces were first proved in \cite{ADM}.)  In the present paper, we carry on to investigate when a mixed Tsirelson space contains higher order $\ell^1$-spreading models {\em hereditarily}.  Again, the first result of this kind is found in \cite{ADM}.  We prove some general characterizations and obtain the result in \cite{ADM} as a corollary.  Roughly speaking, our results show that the complexity of the hereditary finite dimensional $\ell^1$-structure of a mixed Tsirelson space is the same whether it is measured by the existence of higher order $\ell^1$-spreading models or Bourgain's $\ell^1$-index.  These are also related to what may be called ``subsequential minimality" of the mixed Tsirelson space in question and imply that it is arbitrarily distortable.

Denote by $\N$ the set of natural numbers.  For any infinite subset $M$ of $\N$, let $[M]$, respectively $[M]^{<\infty}$, be the set of all infinite and finite subsets of $M$ respectively.  These are subspaces of the power set of $\N$, which is identified with $2^{\N}$ and endowed with the topology of pointwise convergence.  A subset $\F$ of $[\N]^{<\infty}$ is said to be {\em hereditary} if $G \in \F$ whenever $G \subseteq F$ and $F \in \F$.  It is {\em spreading} if for all strictly increasing sequences $(m_i)^k_{i=1}$ and $(n_i)^k_{i=1}$, $(n_i)^k_{i=1} \in \F$ if 
$(m_i)^k_{i=1} \in \F$ and $m_i \leq n_i$ for all $i$.  We also call $(n_i)^k_{i=1}$ a spreading of $(m_i)^k_{i=1}$. A {\em regular} family is a subset of $[\N]^{<\infty}$ that is hereditary, spreading and compact (as a subspace of $2^\N$).  
If $I$ and $J$ are nonempty finite subsets of $\N$, we write $I < J$ to mean $\max I < \min J$.  We also allow that $\emptyset < I$ and $I < \emptyset$.  For a singleton $\{n\}$, $\{n\} < J$ is abbreviated to $n < J$.
Given a family (regular or otherwise) $\F \subseteq [\N]^{<\infty}$, a sequence of sets $(E_i)^k_{i=1}$ is said to be $\F$-{\em admissible} if $(\min E_i)^k_{i=1} \in \F$.  If $\G$ is another family of sets, let 
\[ \F[\G] = \{\cup^k_{i=1}G_i: G_i \in \G, (G_i)^k_{i=1} \text{ is  $\F$-admissible}\} \]
and
\[ (\F,\G) = \{F\cup G: F < G, F \in \F, G \in \G\}.\]
Inductively, set $(\F)^1 = \F$ and $(\F)^{n+1} = (\F,(\F)^n)$ for all $n \in \N$.
It is clear that $\F[\G]$ and $(\F,\G)$ are regular if both $\F$ and $\G$ are.  
A class of regular families that has played a central role is the class of generalized Schreier families \cite{AA}.  The reason for their usefulness as a measure of the complexity of subsets of $[\N]^{<\infty}$ is by now well explained \cite{G, J}.
Let $\cS_0$ consist of all singleton subsets of $\N$ together with the empty set.  Then define $\cS_1$ to be the collection of all $A \in [\N]^{<\infty}$ such that $|A| \leq \min A$ together with the empty set, where $|A|$ denotes the cardinality of the set $A$.  
If $\cS_\alpha$ has been defined for some countable ordinal $\alpha$, set $\cS_{\alpha+1} = \cS_1[\cS_\alpha]$.  For a countable limit ordinal $\alpha$, specify a sequence $(\alpha_n)$ that strictly increases to $\alpha$.  Then define
\[ \cS_\alpha = \{F : F \in \cS_{\alpha_n} \text{ for some $n \leq \min F\} \cup \{\emptyset\}$}. \]
Given a nonempty compact family $\F \subseteq [\N]^{<\infty}$, let $\F^{(0)} = \F$ and $\F^{(1)}$ be the set of all limit points of $\F$.  Continue inductively to derive $\F^{(\alpha+1)} = (\F^{(\alpha)})^{(1)}$ for all ordinals $\alpha$ and $\F^{(\alpha)} = \cap_{\beta<\alpha}\F^{(\beta)}$ for all limit ordinals $\alpha$.  The {\em index} $\io(\F)$ is taken to be the smallest $\alpha$ such that $\F^{(\alpha+1)} = \emptyset$.  
Since $[\N]^{<\infty}$ is countable, $\io(\F) < \omega_1$ for any compact family $\F \subseteq [\N]^{<\infty}$.  
It is well known that $\io(\cS_\alpha) = \omega^\alpha$ for all $\alpha < \omega_1$ \cite[Proposition 4.10]{AA}.

Denote by $c_{00}$ the space of all finitely supported real sequences.  For a finite subset $E$ of $\N$ and $x \in c_{00}$, let $Ex$ be the coordinatewise product of $x$ with the characteristic function of $E$. The sup norm and the $\ell^1$-norm on $c_{00}$ are denoted by $\|\cdot\|_{c_0}$ and $\|\cdot\|_{\ell^1}$ respectively.
Given a sequence $(\F_n)$ of regular families and a nonincreasing null sequence $(\theta_n)^\infty_{n=1}$ in $(0,1)$, 
define a sequence of norms $\|\cdot\|_m$ on $c_{00}$ as follows.  Let
$\|x\|_0 = \|x\|_{c_0}$ and 
\begin{equation}\label{eq 0}
 \|x\|_{m+1} = \max\{\|x\|_m, \sup_n\theta_n\sup\sum^r_{i=1}\|E_ix\|_m\},
\end{equation}
where the last sup is taken over all $\F_n$-admissible sequences $(E_i)^r_{i=1}$.
Since these norms are all dominated by the $\ell^1$-norm, $\|x\| = \lim_m\|x\|_m$ exists and is a norm on $c_{00}$.  The {\em mixed Tsirelson space} $T[(\theta_n,\F_n)^\infty_{n=1}]$ is the completion of $c_{00}$ with respect to the norm $\|\cdot\|$.
From equation (\ref{eq 0}) we can deduce that the norm in $T[(\theta_n,\F_n)^\infty_{n=1}]$ satisfies the implicit equation
\begin{equation}\label{equation 0.1}
\|x\|= \max\{\|x\|_{c_0}, \sup_n\theta_n\sup\sum^r_{i=1}\|E_ix\|\}, 
\end{equation}
with the last sup taken over all $\F_n$-admissible sequences $(E_i)^r_{i=1}$.
{\em For the rest of the paper, we consider a fixed sequence $(\theta_n,\F_n)^\infty_{n=1}$ as above and let $X = T[(\theta_n,\F_n)^\infty_{n=1}]$}.  Set $\alpha_n = \io(\F_n)$ for all $n$.  Families $\F_n$ with $\io(\F_n) = 1$ contains singletons and the empty set only and may be removed without effect on the norm $\|\cdot\|$.  Also the spaces $T[(\theta_n,\F_n)^\infty_{n=1}]$ and $T[(\theta_n,\cup^n_{k=1}\F_k)^\infty_{n=1}]$ are identical (since $(\theta_n)$ is nonincreasing).  
{\em Hence there is no loss of generality in assuming that $\alpha_n > 1$ for all $n$ and that $(\alpha_n)$ is nondecreasing.
We will also assume that $\alpha_n < \sup_m\alpha_m = \omega^{\omega^\xi}$, $0 < \xi < \omega_1$}.  Otherwise, the relevant result has been obtained in \cite[Proposition 2]{LT2}, except for the case when $\xi = 0$.  The coordinate unit vectors $(e_k)$ form an unconditional basis of $X$.  

Given a Banach space $B$ with a basis $(b_k)$, the {\em support} of a vector $x = \sum a_kb_k$ (with respect to $(b_k)$), denoted $\supp x$, is the set of all $k$ such that $a_k \neq 0$.  A {\em block sequence} in $B$ is a sequence $(x_k)$ so that $\supp x_k < \supp x_{k+1}$ for all $k$.  The closed linear span of a block sequence is called a {\em block subspace}.

\section{Technical preliminaries}

In this section, we present some technical results prior to the main discussion.  If $(x_k)$ and $(y_k)$ are sequences of vectors residing in (possibly different) normed spaces, we say that $(x_k)$ {\em dominates} $(y_k)$ if there is a finite positive constant $K$ so that 
\[ \|\sum a_ky_k\| \leq K\|\sum a_kx_k\| \]
for all $(a_k) \in c_{00}$.  Two sequences are {\em equivalent} if they dominate each other.
The first lemma shows that under certain mild assumptions on the families $(\F_n)$, any subsequence of $(e_k)$ is equivalent to its left shift.  The proof uses essentially the idea in \cite[Lemma 2]{CJT}, dressed up in the present language.
The family of all subsets of $\N$ with at most $k$ elements is denoted by $\A_k$.

\begin{lemma} \label{lemma 1}
Assume that for all $n$,  either $\F_n = \A_j$ for some $j \in \N$ or $\F_n[\A_3] \subseteq (\F_n)^2$.
Suppose that $(i_k) \in [\N]$. Let $x = \sum a_ke_{i_{k+1}}$ and $y = \sum a_ke_{i_k}$ for some $(a_k) \in c_{00}$.  Then for any $m$, there exist $E_1 < E_2 < E_3$ such that 
\[ \|x\|_m \leq \sum^3_{i=1}\|E_iy\|_m.\]
Consequently, the sequences $(e_{i_k})$ and $(e_{i_{k+1}})$ are equivalent.
\end{lemma}

\begin{proof}
For any set $E \subseteq \N$, let the left shift of $E$ be the set $L_E = \{i_k: i_{k+1} \in E\}$.  We prove the lemma by induction on $m$.  The case $m = 0$ is clear.  Assume that the lemma holds for some $m$.  If $\|x\|_{m+1} = \|x\|_m$, there is nothing to prove.  Otherwise, $\|x\|_{m+1} = \theta_n\sum^r_{i=1}\|F_ix\|_m$ for some $n$ and some $\F_n$-admissible sequence $(F_i)^r_{i=1}$.  By the inductive hypothesis, there exist $F^i_1 < F^i_2 < F^i_3$ such that 
\[ \|F_ix\|_m \leq \sum^3_{k=1}\|F^i_ky\|_m, \quad 1 \leq i \leq r.\]
We may assume that $F^i_1 \cup F^i_2 \cup F^i_3 \subseteq L_{F_i}$.  We claim that 
$(F^i_k)^{r\ \ \ \ 3}_{i=1\ k=1}$ is $(\A_1\cup \F_n,(\F_n)^2)$-admissible.  
Indeed, if $\F_n = \A_j$ for some $j$, then
\[ (\min F^i_k)^{r\ \ \ \ 3}_{i=1\ k=1} \in \A_j[\A_3] = \A_{3j} = (\F_n)^3 \subseteq (\A_1\cup \F_n,(\F_n)^2).\]
Otherwise,
since $\min F^i_2 \geq \min F_i$, 
\[ \cup^r_{i=1}\{\min F^i_2, \min F^i_3, \min F^{i+1}_1\} \in \F_n[\A_3] \subseteq (\F_n)^2.\]
Clearly, $\{\min F^1_1\} \in \A_1$.  Thus
\[ \cup^r_{i=1}\{\min F^i_1, \min F^i_2, \min F^i_3\} \in  (\A_1\cup \F_n,(\F_n)^2), \]
as claimed.  It follows from the claim that there exist $E_1 < E_2 < E_3$ so that 
$\cup^3_{p=1}E_p = \cup^r_{i=1}\cup^3_{k=1}F^i_k$, 
each $E_p$ is a union of finitely many $F^i_k$ and that $\E_p = \{F^i_k: F^i_k \subseteq E_p\}$ is $(\A_1\cup \F_n)$-admissible if $p = 1$ and $\F_n$-admissible if $p = 2, 3$.  Notice that $\theta_n\sum_{F^i_k\in\E_p}\|F^i_ky\|_m \leq \|E_py\|_{m+1}$ since $\E_p$ is either $\F_n$-admissible or $\A_1$-admissible.  Hence
\begin{equation*}
\|x\|_{m+1} = \theta_n\sum^r_{i=1}\|F_ix\|_m 
\leq \theta_n\sum^r_{i=1}\sum^3_{k=1}\|F^i_ky\|_m 
\leq \sum^3_{p=1}\|E_py\|_{m+1}.
\end{equation*}
Upon taking the limit as $m \to \infty$, we see that $(e_{i_{k+1}})$ is dominated by $(e_{i_k})$.  Since the reverse domination is clear, the two sequences are equivalent.
\end{proof}

A {\em tree} in a Banach space $B$ is a subset $\T$ of $\cup^\infty_{n=1}B^n$ so that $(x_1,\dots,x_n) \in \T$ whenever $(x_1,\dots,x_n,x_{n+1}) \in \T$.  Elements of the tree are called {\em nodes}. It is {\em well-founded} if there is no infinite sequence $(x_n)$ so that $(x_1,\dots,x_m) \in \T$ for all $m$.  If $B$ has a basis, then a tree $\T$ is said to be a {\em block tree} (with respect to the basis) if every node is a block sequence.  For any well-founded tree $\T$, its {\em derived tree} is the tree $\D^{(1)}(\T)$ consisting of all nodes $(x_1,\dots,x_n)$ so that $(x_1,\dots,x_n,x) \in \T$ for some $x$.  Inductively, set $\D^{(\alpha+1)}(\T) = \D^{(1)}(\D^{(\alpha)}(\T))$ for all ordinals $\alpha$ and $\D^{(\alpha)}(\T) = \cap_{\beta<\alpha}\D^{(\beta)}(\T)$ for all limit ordinals $\alpha$.  The {\em order} of a tree $\T$ is the smallest ordinal $o(\T) = \alpha$ such that $\D^{(\alpha)}(\T) = \emptyset$.

\begin{lemma} \label{lemma 2}
Let $\T$ be a well-founded block tree in a Banach space $B$ with a basis.
Define 
\[ \cH = \{(\max\supp x_j)^r_{j=1} : (x_j)^r_{j=1} \in {\T}\} \]
and 
\[ \G = \{G : \text{$G$ is a spreading of a subset of some $H \in \cH$}\}. \]
Then $\G$ is hereditary and spreading.  If $\G$ is compact, then  $\io(\G) \geq o({\T})$.
\end{lemma}

\begin{proof}
It is clear that $\mathcal{G}$ is hereditary and spreading.  Assume that $\G$ is compact. We show by induction on $\xi$ that for all countable ordinal $\xi$, $\io(\G) \geq \xi$ if $o(\T) \geq \xi$. There is nothing to prove if $\xi = 0$.
Suppose the proposition holds for some $\xi<\omega_{1}.$ Let $\T$ be a well-founded 
block tree with $o(\T) \geq \xi+1$.  For each $(x) \in \T,$ let
\[
\T_{x}=\cup_{n=1}^{\infty}\{(x_{1},...,x_{n}) :(x,x_{1},...,x_{n})  \in \T\}.
\]
According to \cite[Proposition 4]{B}, $o(\T) = \sup_{(x)\in \T}(o(\T_{x})+1)$. Therefore,
there exists $(x_{0})  \in \T$ such that $o(\T_{x_{0}}) \geq \xi$. 
By the inductive hypothesis, $\iota(\mathcal{G}') \geq \xi$, where $\G'$ is defined analogously to $\G$ for the tree $\T_{x_0}$. Let $k_{0}=\max\supp x_{0}$. Then
$\{k_{0}\}\cup G\in\mathcal{G}$ whenever $G \in \G'$.
Thus $\{k_{0}\} \in \mathcal{G}^{(\xi)}$.
Since $\mathcal{G}^{(\xi)}$
is spreading, $\{k\} \in \mathcal{G}^{(\xi)}$ for all $k\geq k_{0}$. It follows that
$\iota(\mathcal{G}) \geq \xi+1$.

Suppose $o(\T) \geq \xi_{0},$ where $\xi_{0}$ is a countable limit
ordinal and the proposition holds for all $\xi<\xi_{0}$. Since $o(\T) \geq \xi$ for all $\xi<\xi_{0}$, by the inductive hypothesis,
$\iota(\mathcal{G}) \geq \xi$ for all $\xi<\xi_{0}$. 
Hence $\iota(\mathcal{G})  \geq\xi_{0}$. 
This completes the induction.
\end{proof}

\section{Main results and proofs}

The main results concern two measures of the finite dimensional $\ell^1$-complexity of the space $X$. These are the Bourgain $\ell^1$-index and the existence of $\ell^1$-spreading models of higher order.  Given a finite constant $K$ bigger than $1$, an $\ell^1$-$K$-tree in a Banach space $B$ is a tree in $B$ so that every node $(x_1,\dots,x_n)$ is a normalized sequence such that $\|\sum a_kx_k\| \geq K^{-1}\sum|a_k|$ for all $(a_k)$.  If $B$ has a basis, an $\ell^1$-$K$-block tree is a block tree that is also an $\ell^1$-$K$-tree.  Suppose that $B$ does not contain $\ell^1$, let $I(B,K) = \sup o(\T)$, where the sup is taken over the set of all $\ell^1$-$K$-trees in $X$.  The {\em Bourgain} $\ell^1$-{\em index} is defined to be $I(B) = \sup_{K<\infty}I(B,K)$.  The {\em block indices} $I_b(B,K)$ and $I_b(B)$ are defined analogously using $\ell^1$-block trees.  We refer to \cite{AJO,JO} for thorough investigations of these indices.  In particular, it is shown in \cite{JO} that for a Banach space $B$ with a basis, $I_b(B) = I(B)$ if either one is $\geq \omega^\omega$.  With the same notation as above, a normalized sequence $(x_k)$ is said to be an $\ell^1$-$\cS_{\beta}$-spreading model with constant $K$ if $\|\sum_{k\in F}a_kx_k\| \geq K^{-1}\sum_{k\in F}|a_k|$ whenever $F \in \cS_\beta$.

We are now ready to work our way towards the main Theorem \ref{theorem 8}.  The major parts of the computations are contained in Proposition \ref{proposition 5} and Lemma \ref{lemma 8} (tree splitting lemma).
Let $(y_k)$ be a normalized block sequence in $X$ and let $Y$ be the block subspace $[(y_k)]$.  For any $n \in \N$, we call the space $Y_n = [(y_k)^\infty_{k=n}]$ the $n$-tail of $Y$. We emphasize that in the next lemma both admissibility and the support of a vector are taken with respect to the basis $(e_k)$.  Recall the assumption that $(\alpha_n) = (\io(\F_n))$ is a nondecreasing sequence which converges to $\omega^{\omega^\xi}$ nontrivially.

\begin{lemma} \label{lemma 3}
Assume that $I_b(Y) > \omega^{\omega^\xi}$.  Then there exists a constant $C < \infty$ such that for all $n \in \N$, there exists a normalized vector $x$ in the $n$-tail of $Y$ such that $\sum\|E_ix\| \leq C$ whenever $(E_i)$ is $\F_k$-admissible for some $k \leq n$.
\end{lemma}

\begin{proof}
There exists $K < \infty$ such that $I_b(Y,K) \geq \omega^{\omega^\xi}$.  Let $\T$ be an $\ell^1$-$K$-block tree in $Y$ such that $o(\T) \geq \omega^{\omega^\xi}$.  Given $n$, consider the tree $\widehat{\T}$ consisting of all nodes of the form $(x_j)^r_{j=n}$ for some $(x_j)^r_{j=1} \in \T$, $r \geq n$.  Then $\widehat{\T}$ is an $\ell^1$-$K$-block tree in $Y_n$ such that $o(\widehat{\T}) \geq \omega^{\omega^\xi}$. 
Choose $\alpha$ and $\beta$ so that $\alpha_n < \omega^\alpha < \omega^\beta < \omega^{\omega^\xi}$.  Define 
\[ \cH = \{(\max\supp x_j)^r_{j=n} : (x_j)^r_{j=n} \in \widehat{\T}\} \]
and 
\[ \G = \{G : \text{$G$ is a spreading of a subset of some $H \in \cH$}\}. \]
By Lemma \ref{lemma 2}, $\G$ is hereditary and spreading, and either $\G$ is noncompact or it is compact with $\io(\G) \geq o(\widehat{\T}) \geq {\omega^{\omega^\xi}}$.  By \cite[Theorem 1.1]{G}, there exists $M \in [\N]$ such that
\[ \cup^n_{k=1}\F_k \cap [M]^{<\infty} \subseteq \cS_\alpha \cap [M]^{<\infty}
\subseteq \cS_\beta \cap [M]^{<\infty} \subseteq \G.\]
Now \cite[Proposition 3.6]{OTW} gives a finite set $G \in \cS_\beta \cap [M]^{<\infty}$ and a sequence of positive numbers $(a_p)_{p\in G}$ such that $\sum a_p = 1$ and $\sum_{p\in F}a_p < \theta_n$ whenever $F \subseteq G$ and $F \in \cS_\alpha$.
By definition, there exist a node $(x_j)^r_{j=n} \in \widehat{\T}$ and a subset $J$ of the integer interval $[n,r]$ such that $G$ is a spreading of $(\max\supp x_j)_{j\in J}$.
Denote the unique order preserving bijection from $J$ onto $G$ by $u$ and consider the vector $y = \sum_{j\in J}a_{u(j)}x_j$.  Since $(x_j)^r_{j=n}$ is a normalized $\ell^1$-$K$-block sequence in $Y_n$ and $\sum a_{u(j)} = 1$, $y \in Y_n$ and $\|y\| \geq 1/K$.
Let $(E_i)$ be $\F_k$-admissible for some $k \leq n$.  For each $j\in J$, let $\E_j$ be the collection of all $E_i$'s that have nonempty intersection with $\supp x_{j'}$ if and only if $j' = j$. Also let $\E'$ be the collection of all $E_i$ such that $E_i$ intersects $\supp x_j$ for at least two $j\in J$.  Since $(E_i)$ is $\F_k$-admissible, for each $j \in J$,
\[ \sum_{E_i\in\E_j}\|E_iy\| \leq a_{u(j)}\theta_k^{-1}\|x_j\| = a_{u(j)}\theta_k^{-1}. \]
Set $J' = \{j\in J: \E_j \neq \emptyset\}$.  The $\F_k$-admissiblity of $(E_i)$ implies that
$(\max\supp x_j)_{j\in J'} \in \F_k$.  Thus $u(J')$, being a spreading of this set, also belongs to $\F_k$.  Since $u(J') \subseteq G \in [M]^{<\infty}$, we conclude that $u(J') \in \F_k \cap [M]^{<\infty} \subseteq \cS_\alpha$.  Hence $\sum_{j\in J'}a_{u(j)} < \theta_n$.
Also, since each $\supp x_j$, $j \in J$, intersects at most two $E_i$ in $\E'$,
\[ \sum_{E_i\in\E'}\|E_iy\| \leq \sum_{j\in J}a_{u(j)}\sum_{E_i\in \E'}\|E_ix_j\| \leq 2\sum_{j\in J}a_{u(j)} = 2.\]
Therefore,
\begin{align*}
\sum\|E_iy\| &= \sum_{E_i\in\E'}\|E_iy\| + \sum_{j\in J'}\sum_{E_i\in\E_j}\|E_iy\| \\
& \leq 2 + \theta_k^{-1}\sum_{j\in J'}a_{u(j)} \leq 3.
\end{align*}
It is clear that the normalized element $x = y/\|y\|$ satisfies the statement of the lemma with the constant $C = 3K$.
\end{proof}

We pause to introduce another method of computing the norm of an element in $X$ using norming trees.  This is derived from the implicit description  of the norm in $X$ (equation (\ref{equation 0.1})) and have been used in \cite{Be, LT2, OT}.   
An ($(\F_k)$-){\em admissible tree} is a finite collection of elements $(E_{i}^{m})$,
$0\leq m\leq r,$ $1\leq i\leq k(m)$, in $[\N]^{<\infty}$ with the following properties.
\begin{enumerate}
\item $k(0) = 1$, 
\item For each $m$, $E_{1}^{m}<E_{2}^{m}<\dots<E_{k(m)}^{m}$,
\item Every $E_{i}^{m+1}$ is a subset of some $E_{j}^{m}$,
\item For each $j$ and $m$, the collection $\{E_i^{m+1}: E_i^{m+1} \subseteq E^m_j\}$ is $\F_k$-admissible for some $k$.
\end{enumerate}
The set $E^0_1$ is called the {\em root} of the admissible tree.
The elements $E_{i}^{m}$ are called {\em nodes} of the
tree. If $E_{i}^{n}\subseteq E_{j}^{m}$ and $n>m$, we
say that $E_{i}^{n}$ is a {\em descendant} of $E_{j}^{m}$ and $E_{j}^{m}$ is
an {\em ancestor} of $E_{i}^{n}$. If, in the above notation, $n=m+1$, then
$E_{i}^{n}$ is said to be an {\em immediate successor} of $E_{j}^{m}$, and
$E_{j}^{m}$ the {\em immediate predecessor} of $E_{i}^{n}$. Nodes with no
descendants are called {\em terminal nodes} or {\em leaves} of the tree.
Assign {\em tags} to the individual nodes inductively as follows. Let 
$t(E_{1}^{0}) = 1$.  If $t(E_{i}^{m})$ has
been defined and the collection $(E_{j}^{m+1})$ of all
immediate successors of $E_{i}^{m}$ forms an $\F_k$-admissible
collection, then define $t(E_{j}^{m+1}) = \theta_kt(E^m_i)$ for all immediate successors $E_{j}^{m+1}$ of $E_{i}^{m}.$  If $x \in c_{00}$ and $\T$ is an admissible tree, let $\T x = \sum t(E)\|Ex\|_{c_0}$ where the sum is taken over all leaves in
$\T$. It follows from the implicit description (equation (\ref{equation 0.1})) of the norm in $X$ that $\|x\| = \max \T x$, with the maximum taken over the set of all admissible trees.
Let us also point out that if $\E$ is a collection of pairwise disjoint nodes of an admissible tree $\T$ so that $E \subseteq \cup \E$ for every leaf $E$ of $\T$ and $x \in c_{00}$, then $\T x = \sum_{F \in \E}t(F)\|Fx\|$.

Suppose that $\F$ and $\G$ are families of finite subsets of $\N$.  An element $G \in \G$ is {\em maximal} (in $\G$) if it is not properly contained in any other element in $\G$.  Define $\F\ominus \G$ to be the collection of all sets $F$ so that there is a maximal $G \in \G$, $G < F$, with $G\cup F \in \F$.  We say that a sequence of regular families $(\F_n)$ is {\em tame} if 
(1) for each $n$, either $\F_n = \A_j$ for some $j$ or $\F_n[\A_3] \subseteq (\F_n)^2$ and
(2) there exists $n_0 \in \N$ so that $(\F_n\ominus\F_{n_0})[\A_2] \subseteq  \F_n$ whenever $n > n_0$.

\begin{proposition}\label{proposition 5}
Assume that $(\F_n)$ is a tame sequence.  Let $Y$ be a block subspace of $X$.  Suppose that there exists a constant $C < \infty$ such that for all $n \in \N$, there is a normalized vector $x$ in the $n$-tail of $Y$ such that $\sum\|E_ix\| \leq C$ whenever $(E_i)$ is $\F_k$-admissible for some $k \leq n$.  Then there exists a normalized block sequence $(z_n)$ in $Y$ that is equivalent to a subsequence of $(e_k)$.
\end{proposition}

\begin{proof}
Let $n_0$ be the integer occurring in the definition of tameness for the sequence $(\F_n)$.  Inductively, choose a normalized block sequence $(z_n)$ in $Y$ and a strictly increasing sequence $(m_n)^\infty_{n=0}$ in $\N$ so that $m_0 > n_0$, $\theta_{m_n}\|z_n\|_{\ell^1} \leq 2^{-n}$ and $\sum\|E_iz_n\| \leq C$ whenever $(E_i)$ is $\cup^{m_{n-1}}_{r=1}\F_{r}$-admissible, $n \in \N$.  Consider $z = \sum a_nz_n$ for some $(a_n) \in c_{00}$ and let $y = \sum a_ne_{k_n}$, where $k_n = \max\supp z_n$.  Let $\T$ be an admissible tree that norms $z$.  Without loss of generality, we may assume that all nodes in $\T$ are integer intervals and that all leaves in $\T$ are singletons.  Say that a node is short if it intersects $\supp z_n$ for exactly one $n$.  On the other hand, call a node long if it intersects $\supp z_n$ for more than one $n$.  
The tree $\T$ is endowed with the natural partial order of reverse inclusion.
Let $\E$ be the collection of all minimal short nodes in $\T$.  Then $\|z\| = \sum_{E\in\E}t(E)\|Ez\|$.  For each $n$, let $\E_n$ be the collection of all nodes in $\E$ that intersects only $\supp z_n$.  In particular, $\E = \cup \E_n$.  Further subdivide each set $\E_n$ into two subsets $\E_n'$ and $\E_n''$ depending on whether $t(E) \leq \theta_{m_n}$ or not.  We have
\begin{equation}\label{eq 1}
\sum_n\sum_{E\in\E'_n}t(E)\|Ez\| \leq \sum_n\theta_{m_n}|a_n|\|z_n\|_{\ell^1} \leq \sum_n\frac{|a_n|}{2^n} \leq \|y\|.
\end{equation}
For each $n$, let $\D_n$ be the set of all minimal elements in the set of all nodes in $\T$ that are immediate predecessors of some node in $\E''_n$.  Since $\D_n$ consists of pairwise disjoint long nodes that intersect $\supp z_n$, $|\D_n| \leq 2$ for all $n$.  For each $D \in \D_n$, let $\E''_n(D) = \{E \in \E''_n: E \subseteq D\}$ and let $\tilde{\E}_n''(D)$ be the subset of $\E''_n$ consisting of all $E \in \E''_n$ that are immediate successors of $D$.  Fix $E_{n,D} \in \tilde{\E}_n''(D)$ and $j_{n,D} \in E_{n,D}\cap \supp z_n$ arbitrarily and set $w = \sum_n\sum_{D\in\D_n}a_ne_{j_{n,D}}$.  Since $|\D_n| \leq 2$ for all $n$, $\|w\| \leq 2\|y\|$.  Any immediate successor of $D$ that contains some $E \in \E''_n(D)\backslash\tilde{\E}''_n(D)$ must be a long node.  Hence there are at most two immediate successors of $D$, say $G_1$ and $G_2$, that all nodes in 
$\E''_n(D)\backslash\tilde{\E}''_n(D)$ are descended from.  Note that $t(G_1) = t(G_2) = t(E_{n,D})$ since they are all immediate successors of the same node.  Thus
\begin{equation*}
\sum_{E\in \E''_n(D)\backslash\tilde{\E}''_n(D)}t(E)\|Ez_n\| \leq \sum^2_{i=1}t(G_i)\|G_iz_n\| 
\leq 2t(E_{n,D}).
\end{equation*} 
Hence
\begin{align}\label{eq 2}
\sum_n\sum_{D\in\D_n}&\sum_{E\in \E''_n(D)\backslash\tilde{\E}''_n(D)}t(E)\|Ez\|\\ \notag 
&= \sum_n\sum_{D\in\D_n}\sum_{E\in \E''_n(D)\backslash\tilde{\E}''_n(D)}t(E)|a_n|\|Ez_n\|\\ \notag
&\leq \sum_n\sum_{D\in\D_n}2|a_n|t(E_{n,D})\\ \notag
&\leq 2\sum_n\sum_{D\in\D_n}t(E_{n,D})\|E_{n,D}w\| \leq 2\|w\| \leq 4\|y\|.
\end{align}
Now let $\D'_n$ consists of those $D$ in $\D_n$ such that $\tilde{\E}''_n(D)$ is $\cup^{m_{n-1}}_{r=1}\F_{r}$-admissible.  Then
\begin{align}\label{eq 3}
\sum_n\sum_{D\in\D'_n}&\sum_{E\in\tilde{\E}''_n(D)}t(E)\|Ez\|\\ \notag
&= \sum_n\sum_{D\in\D'_n}\sum_{E\in\tilde{\E}''_n(D)}t(E)|a_n|\|Ez_n\|\\ \notag
&\leq C\sum_n\sum_{D\in\D'_n}t(E_{n,D})|a_n| \\ \notag
&\leq C\sum_n\sum_{D\in\D'_n}t(E_{n,D})\|E_{n,D}w\| \leq C\|w\| \leq 2C\|y\|.
\end{align}
It remains to consider the nodes that belong to $\D_n\backslash \D'_n$ for some $n$. We have
\begin{align*}\label{eq 4}
\sum_n\sum_{D\in\D_n\backslash\D'_n}&\sum_{E\in\tilde{\E}''_n(D)}t(E)\|Ez\| \\ \notag
&=  \sum_n\sum_{D\in\D_n\backslash\D'_n}\sum_{E\in\tilde{\E}''_n(D)}t(E)|a_n|\|Ez_n\| \\ \notag
&\leq \sum_n\sum_{D\in\D_n\backslash\D'_n}t(D)|a_n|\|Dz_n\|\\ \notag
&\leq \sum_n\sum_{D\in\D_n\backslash\D'_n}t(D)|a_n|. 
\end{align*}
But by Lemma \ref{lemma 8} below, 
\begin{equation*} \label{5}
\sum_n\sum_{D\in\D_n\backslash\D'_n}t(D)|a_n| \leq 4\|y\|.
\end{equation*}
Thus 
\begin{equation}\label{eq 6}
\sum_n\sum_{D\in\D_n\backslash\D'_n}\sum_{E\in\tilde{\E}''_n(D)}t(E)\|Ez\| \leq 4\|y\|.
\end{equation}
Combining inequalities (\ref{eq 1}) to (\ref{eq 6}), we see that 
\[ \|z\| = \sum_{E\in\E}t(E)\|Ez\|\leq (9+2C)\|y\|.\]
Hence $(z_n)$ is dominated by $(e_{k_n})$, where $k_n = \max\supp z_n$.  On the other hand, $(z_n)$ dominates $(e_{k_{n-1}})$ (take $k_0 =1$).  Therefore, using the tameness of $(\F_n)$, we see that $(z_n)$ is equivalent to $(e_{k_n})$ by Lemma \ref{lemma 1}.
\end{proof}

\begin{lemma}\label{lemma 7}
Suppose that $n_1 < n_2$ and $D \in \D_{n_2}\backslash\D'_{n_2}$.  Then no descendant of $D$ belongs to ${\E}''_{n_1}$.  In particular, $D\notin \D_{n_1}$.
\end{lemma}

\begin{proof}
If $E$ is a descendant of $D \in \D_{n_2}\backslash\D'_{n_2}$, then $t(E) \leq t(F)$ for any immediate successor $F$ of $D$.  In particular, $t(E) \leq t(F)$ for all $F \in \tilde{\E}''_{n_2}(D)$.  By definition of $\D'_{n_2}$, $\tilde{\E}''_{n_2}(D)$ is not $\F_r$-admissible for all $r \leq {m_{n_1}}$.  Hence $t(F) < \theta_{m_{n_1}}$ for all $F \in \tilde{\E}''_{n_2}(D)$.  Therefore, $t(E) < \theta_{m_{n_1}}$ if $E$ is a descendant of $D \in \D_{n_2}\backslash\D'_{n_2}$.  This shows that $E \notin 
{\E}''_{n_1}$ by definition of $\E''_{n_1}$.
\end{proof}

Let $\T'$ be the subtree of $\T$ consisting of all nodes in $\tilde{\D} = \cup_n(\D_n\backslash\D'_n)$ and their ancestors.  By Lemma \ref{lemma 7}, for each $D \in \tilde{\D}$, there is a unique $n = n_D$ such that $D \in \D_n\backslash\D'_n$. If $G$ is a node in $\T'$, let $\tilde{\D}(G)$ consists of all $D \in \tilde{\D}$ such that $D \subseteq G$.  Recall the vector $w$ defined in the proof of Proposition \ref{proposition 5} above.  It was observed that $\|w\| \leq 2\|y\|$.

\begin{lemma}\label{lemma 8}
For any $G \in \T'$, there exist subsets $G_1$ and $G_2$ of $G$, $G_1 < G_2$, and admissible trees $\T_1$ and $\T_2$ with roots $G_1$ and $G_2$ respectively so that
\[ \sum_{D\in\tilde{\D}(G)}t(D)|a_{n_D}| \leq t(G)(\T_1w+\T_2w).\]
In particular, 
\[ \sum_n\sum_{D\in\D_n\backslash\D'_n}t(D)|a_n| \leq 4\|y\|. \]
\end{lemma}

\begin{proof}
The second inequality follows from the first by taking $G$ to be the root of $\T'$ (which is also the root of $\T$).  To prove the first inequality, we begin at the terminal nodes of $\T'$ and work our way up the tree.  Let $G$ be a terminal node of $\T'$. Then $G \in \tilde{\D}$.  In this case, take $G_1 = [1,\max\supp z_{n_G}]\cap G$ and $G_2 = G\backslash G_1$.  Clearly, $G_1$ and $G_2$ are subsets of $G$ such that $G_1 < G_2$.  Set $\T_1$ and $\T_2$ to be the trivial trees $\T_i = \{G_i\}$, $i = 1, 2$.  Now
\[ \sum_{D\in\tilde{\D}(G)}t(D)|a_{n_D}| = t(G)|a_{n_G}| \leq t(G)\|G_1w\|_{c_0} = t(G)\T_1w.\]
Thus the lemma holds in this case.

Next, take a node $G \in \T'$ and assume that the lemma has been proved for all descendants of $G$ in $\T'$.  List the immediate successors of $G$ in $\T'$ from left to right as $\{H_1,\dots,H_r\}$.  By the assumption, for each $j$, $1\leq j \leq r$, there are subsets $H^i_j$ of $H_j$, and admissible trees $\T^i_j$, $i = 1, 2$, such that $H^1_j < H^2_j$, the root of $\T^i_j$ is $H^i_j$ and 
\[ \sum_{D\in\tilde{\D}(H_j)}t(D)|a_{n_D}| \leq t(H_j)(\T^1_jw+\T^2_jw).\]
We divide the rest of the proof into two cases.\\

\noindent\underline{Case 1}. $G \in \tilde{\D}$.

\noindent The sets in the collection $\tilde{\E}''_{n_G}(G) \cup \{H_j\}^r_{j=1}$  are all immediate successors of $G$ in the tree $\T$.  We claim that $E < H_1$ for any $E \in \tilde{\E}''_{n_G}(G)$.  Indeed either $H_1$ or a descendant of $H_1$ belongs to $\tilde{\D}$.  Denote this node by $I$.  Thus $G \in \D_{n_G}\backslash\D'_{n_G}$ has a descendant in ${\E}''_{n_I}$.  By Lemma \ref{lemma 7}, $n_I \geq n_G$.  Since $I \subsetneq G$, $n_I \neq n_G$ by the minimality condition in the definition of $\D_n$.  Hence $n_I > n_G$.  Now any $E$ in $\tilde{\E}''_{n_G}(G)$ intersects only $\supp z_{n_G}$ while $H_1$ must intersect $\supp z_{n_I}$.  Therefore, $E < H_1$, as claimed. To continue with the proof, set $G_1 = G\cap [1,k]$, where $k = \max\cup\tilde{\E}''_{n_G}(G)$, and $G_2 = G\backslash G_1$.  Then take $\T_1$ to be the trivial tree $\{G_1\}$ and $\T_2$ to be the tree $\{G_2\}\cup (\cup_{i,j}\T^i_j)$.  The admissibility of $\T_1$ is clear.  To verify the admissibility of $\T_2$, it suffices to show the admissibility of the decomposition of $G_2$ into $\{H^i_j\}_{i,j}$.  Since $\tilde{\E}''_{n_G}(G) \cup \{H_j\}^r_{j=1}$  are all immediate successors of $G$ in the tree $\T$, the collection is $\F_n$-admissible for some $n$. However, 
$\tilde{\E}''_{n_G}(G)$ is not $\F_{r}$-admissible for any $r \leq m_{n_G-1}$. Thus $n > m_{n_G-1} > n_0$ and $(\min H_j) \in \F_n \ominus \F_{n_0}$.  By the tameness of $(\F_n)$, $(\min H^i_j) \in (\F_n \ominus \F_{n_0})[\A_2] \subseteq \F_n$.  Hence $(H^i_j)$ is $\F_n$-admissible, as required.  Now
\begin{align*}
\T_1w &= \|G_1w\|_{c_0} \geq |a_{n_G}| \\
\intertext{and}
\T_2w &= \theta_n\sum_{i,j}\T^i_jw =  \sum_{i,j}\frac{t(H_j)}{t(G)}\T^i_jw \geq \sum_j\sum_{D\in\tilde{\D}(H_j)}\frac{t(D)}{t(G)}|a_{n_D}|.
\end{align*}
Therefore,
\begin{align*}
\sum_{D\in\tilde{\D}(G)}t(D)|a_{n_D}| &= t(G)|a_{n_G}| + \sum_j\sum_{D\in\tilde{\D}(H_j)}t(D)|a_{n_D}|\\
&\leq t(G)(\T_1w+\T_2w).
\end{align*}

\noindent\underline{Case 2}. $G \notin \tilde{\D}$.\\

\noindent Suppose that in the tree $\T$, the immediate successors of $G$ form an $\F_n$-admissible collection.  In particular, $\{H_j\}^r_{j=1}$ is $\F_n$-admissible.  
We claim that $(\min H^i_j) \in (\F_n)^2$.  This is clear if $\F_n = \A_j$ for some $j$.  Otherwise, $(\min H^i_j) \in \F_n[\A_2] \subseteq (\F_n)^2$ by the tameness of $(\F_n)$. Choose index sets $I_1$ and $I_2$ such that 
$I_1\cup I_2 = \{(i,j): 1\leq i\leq 2, 1 \leq j\leq r\}$,
$\{H^i_j: (i,j)\in I_k\}$ is $\F_n$-admissible, $k = 1,2$, and that $H^i_j < H^{i'}_{j'}$ whenever $(i,j) \in I_1$ and $(i',j')\in I_2$.  Set $G_1 = G \cap [1,p]$, where $p = \max\cup\{H^i_j: (i,j)\in I_1\}$ and $G_2 = G\backslash G_1$.  Define $\T_k$ to be the tree $\{G_k\} \cup (\cup_{(i,j)\in I_k}\T^i_j)$, $k = 1, 2$.  The admissibility of $\T_1$ and $\T_2$ follows by construction.  Finally,
\begin{align*}
t(G)\sum_k\T_kw &= t(G)\theta_n\sum_k\sum_{(i,j)\in I_k}\T^i_jw \\
&= \sum_{i,j}{t(H_j)}\T^i_jw \\
&\geq \sum_{j}\sum_{D\in\tilde{\D}(H_j)}{t(D)}|a_{n_D}|\\
& = \sum_{D\in\tilde{\D}(G)}{t(D)}|a_{n_D}|.
\end{align*}
\end{proof}

Given a nonzero ordinal $\alpha$ with Cantor normal form $\omega^{\beta_1}\cdot m_1+\dots+ \omega^{\beta_n}\cdot m_n$, let $\ell(\alpha) = \beta_1$.  For any $m \in \N$ and $\ep > 0$, define
\[ \gamma(\ep,m) = \max\{\ell(\alpha_{n_s}\cdots\alpha_{n_1}): \ep\theta_{n_1}\cdots\theta_{n_s} > \theta_m\}\ (\max \emptyset = 0).\]
The sequence $((\theta_n,\F_n))^\infty_{n=1}$ is said to satisfy $(\dagger)$ if there exists $\ep > 0$ such that for all $\beta< \omega^\xi$, there exists $m \in \N$ such that $\gamma(\ep,m) + 2+ \beta < \ell(\alpha_m)$.

\begin{theorem} \cite[Theorems 4 and 12]{LT2} \label{theorem 7}
Assume that $(\dagger)$ holds.  Then for any $M \in [\N]$, $[(e_k)_{k\in M}]$ contains an $\ell^1$-$\cS_{\omega^\xi}$-spreading model.  On the other hand, if $(\dagger)$ fails, then for all $M \in [\N]$, there exists $N \in [M]$ such that $I_b([(e_k)_{k\in N}]) = \omega^{\omega^\xi}$.
\end{theorem}

Recall that a Banach space $(B,\|\cdot\|)$ is said to be $\lambda$-{\em distortable} if there is an equivalent norm $|\cdot|$ on $B$ so that for every infinite dimensional subspace $Y$ of $X$, there are $\|\cdot\|$-normalized vectors $y$ and $z$ in $Y$ so that $|y|/|z| > \lambda$.  A space is {\em arbitrarily distortable} if it is $\lambda$-distortable for all $\lambda > 1$.

\begin{theorem}\label {theorem 8}
Assume that $(\F_n)$ is a tame sequence.
The following statements are equivalent for any block subspace $Y$ of $X$.  
\begin{enumerate}
\item Property $(\dagger)$ holds and every block subspace $Z$ of $Y$ contains a block sequence equivalent to a subsequence of $(e_k)$.
\item Every block subspace $Z$ of $Y$ contains an $\ell^1$-$\cS_{\omega^\xi}$-spreading model.
\item The Bourgain $\ell^1$-index $I_b(Z) = I(Z) > \omega^{\omega^\xi}$ for any block subspace $Z$ of $Y$.
\end{enumerate}
Moreover, if one (and hence all) of the equivalent conditions holds for a block subspace $Y$ of $X$, then $Y$ is arbitrarily distortable.
\end{theorem}

\begin{proof}
The implication (1) $\implies$ (2) follows from the first part of Theorem \ref{theorem 7}.  Let $Z$ be a block subspace of $Y$. If (2) holds, then $I(Z,K) \geq \omega^{\omega^\xi}$ for some $K < \infty$.  By \cite[Lemma 5.7]{JO}, $I_b(Z) = I(Z) > \omega^{\omega^\xi}$.  Assume that condition (3) holds.  By Lemma \ref{lemma 3} and Proposition \ref{proposition 5}, $Z$ contains a normalized block sequence equivalent to a subsequence of $(e_k)$.  Say $(z_n)$ is a normalized block sequence in $Z$ equivalent to $(e_k)_{k\in M}$ for some $M \in [\N]$.  If $(\dagger)$ fails, by the second part of Theorem \ref{theorem 7}, there exists $N \in [M]$ such that $I_b([(e_k)_{k\in N}]) = \omega^{\omega^\xi}$.  Hence $I_b([(z_{n_j})]) = \omega^{\omega^\xi}$ for some subsequence $(z_{n_j})$ of $(z_n)$.  This contradicts (3) since $[(z_{n_j})]$ is a block subspace of $Y$. This proves condition (1).

Assume that the conditions hold for a block subspace $Y$ of $X$.  For each $n$, consider the equivalent norm $|||\cdot|||_n$ on $X$ defined by 
\[ |||x|||_n = \sup\{\sum\|E_ix\|: (E_i) \text{ is $\F_n$-admissible}\}.\]
Let $Z$ be a block subspace of $Y$.  By condition (3) and Lemma \ref{lemma 3}, there exists $C_1 < \infty$ such that for all $n$, there exists $z \in Z$ such that $\|z\| = 1$ and $|||z|||_n \leq C_1$.
On the other hand, by condition (1), $Z$ contains a normalized block sequence $(z_k)_{k\in M}$ that is $C_2$-equivalent to a subsequence $(e_{k})_{k\in M}$ of $(e_k)$.  
Let $\ep$ be the constant given by property $(\dagger)$.  It follows from property $(\dagger)$ that there are infinitely many $m$ such that $\gamma(\ep,m)+2 < \ell(\alpha_m)$.  Fix such an $m$ and let $\gamma = \gamma(\ep,m)$.  By \cite[Theorem 1.1]{G}, there exists $N \in [M]$ such that $\cS_{\gamma+2} \cap [N]^{<\infty} \subseteq \F_m$.
By \cite[Lemma 19]{LT}, there exists $x \in c_{00}$ such that $\|x\| \leq 1+ 1/\ep$, $\|x\|_{\ell^1} = \theta^{-1}_m$ and $\supp x \in \cS_{\gamma+2}\cap [N]^{<\infty}$.  Say $x = \sum_{k\in I}a_ke_k$ for some $I \in [N]^{<\infty}$.  Consider the corresponding element $y = \sum_{k\in I}a_kz_k/\|\sum_{k\in I}a_kz_k\|$.  Since $(z_k)_{k\in N}$ is $C_2$-equivalent to $(e_k)_{k\in N}$, 
\[ \|\sum_{k\in I}a_kz_k\| \leq C_2\|x\| \leq C_2(1+\ep^{-1}). \] 
For each $k$, let $E_k = \supp z_k$.  Then $(\min E_k)_{k\in I}$ is a spreading of $(e_k)_{k\in I} = \supp x$.  Hence $(E_k)_{k\in I}$ is $\F_m$-admissible since $\supp x \in \F_m$.  Therefore,
\[ |||\sum_{k\in I}a_kz_k|||_m \geq \sum_{k\in I}\|E_k\sum_{j\in I}a_jz_j\| = \sum_{k\in I}|a_k| = \|x\|_{\ell^1} = \theta^{-1}_m.\]
Hence $|||y|||_m \geq C^{-1}_2(1+\ep^{-1})^{-1}\theta^{-1}_m$.  
The existence of $z$ and $y$ shows that $Y$ is $C^{-1}_1C^{-1}_2(1+\ep^{-1})^{-1}\theta^{-1}_m$-distortable.  Since this holds for infinitely many $m$, $Y$ is arbitrarily distortable.
\end{proof}

\begin{corollary}
Assume that  $(\F_n)$ is a tame sequence. 
If $\xi$ is a limit ordinal, the following statements hold.
\begin{enumerate}
\item Every block subspace of $X$ contains an $\ell^1$-$\cS_{\omega^\xi}$-spreading model.
\item Every block subspace of $X$ contains a block sequence equivalent to a subsequence of $(e_k)$.
\item $X$ is arbitrarily distortable.
\end{enumerate}
\end{corollary}

\begin{proof}
If $(z_n)$ is a normalized block sequence in $X$, 
and $F$ is a set such that $\{\min \supp z_n\}_{n\in F}\in \F_m$,
then $\|\sum a_nz_n\| \geq \theta_m\sum_F|a_n|$.  In particular, $I_b(Y,\theta_m^{-1}) \geq {{\alpha_m}}$ for all block subspaces $Y$ of $X$ and all $m$.  By the proof of Theorem 1.1 in \cite{JO}, if $I_b(Y,K) \geq \alpha^2$, then $I_b(Y,\sqrt{K}) \geq \alpha$.  Now for any $\beta < \omega^\xi$, there exists $m$ such that $\omega^{\beta\cdot\omega} < \alpha_m$.  Thus $(\omega^\beta)^{2^k} < \alpha_m$  for all $k$.  It follows that $I_b(Y,\theta^{-1/2^k}_m) \geq \omega^\beta$.  Hence $I_b(Y,1+\ep) \geq \omega^\beta$ for any $\ep > 0$ and any $\beta < \omega^\xi$.  Therefore, $I_b(Y,1+\ep) \geq \omega^{\omega^\xi}$ for any $\ep > 0$.  By \cite[Lemma 5.7]{JO}, $I_b(Y) > \omega^{\omega^\xi}$.  The conclusions of the corollary now follows from Theorem \ref{theorem 8}.
\end{proof}

\begin{proposition}
The sequence $(\cS_{\beta_n})$ is tame for any sequence of nonzero countable ordinals $(\beta_n)$. 
\end{proposition}

\begin{proof}
Let $\alpha$ be a nonzero countable ordinal.  The fact that $\cS_\alpha[\A_3] \subseteq (\cS_\alpha)^2$ was shown in the Remark following Proposition 9 in \cite{LT}.  We show that $(\cS_\alpha\ominus\cS_1)[\A_2] \subseteq \cS_\alpha$ by induction on $\alpha$.  If $\alpha = 1$, this is clear.  Assume that the inclusion holds for some $\alpha$.  Suppose $E \in (\cS_{\alpha+1}\ominus\cS_1)[\A_2]$.  Then $E = \cup^k_{i=1}E_i$, $E_1 < \dots < E_k$, $E_i \in \A_2$, and $F = \{\min E_i\}^k_{i=1} \in \cS_{\alpha+1}\ominus\cS_1$.  
There is a maximal $\cS_1$ set $G$ such that $G < F$ and $G \cup F \in \cS_{\alpha+1}$.  Let $\min G = n$.  Then $|G| = n$ and hence $\min F \geq 2n$.  Note that $F \subseteq G\cup F \in \cS_{\alpha+1}$.  Thus we may 
write $F$ as $\cup^r_{j=1}H_j$, where $H_1 < \dots < H_r$, $H_j \in \cS_\alpha$, and $r \leq n$.  
Since $\cS_\alpha[\A_2] \subseteq (\cS_\alpha)^2$, $\cup\{E_i: \min E_i \in H_j\} \in (\cS_\alpha)^2$ for all $j$.
Therefore, 
\[ E \subseteq \cup^r_{j=1}\cup\{E_i: \min E_i \in H_j\} \in (\cS_\alpha)^{2r} \]
and $2r \leq 2n \leq \min F = \min E$.  Hence $E \in \cS_{\alpha+1}$, as required.

Finally, suppose the inclusion holds for all $\alpha' < \alpha$, where $\alpha$ is a limit ordinal.  Let $(\alpha_n)$ be the sequence of ordinals used to define $\cS_\alpha$.  If $E \in (\cS_\alpha\ominus\cS_1)[\A_2]$, then $E \in (\cS_{\alpha_n}\ominus\cS_1)[\A_2]$ for some $n \leq \min E$.  Thus $E \in \cS_{\alpha_n}$ for some $n \leq \min E$.  Hence $E \in \cS_\alpha$.

Observe that $\cS_1\subseteq \cS_\alpha$ for any nonzero countable ordinal $\alpha$. Therefore, if $n > 1$,
\[ (\cS_{\beta_n}\ominus\cS_{\beta_1})[\A_2] \subseteq (\cS_{\beta_n}\ominus\cS_{1})[\A_2] \subseteq \cS_{\beta_n}.\]
\end{proof}

\begin{theorem}\label{theorem 10}
Let $(\theta_n)$ be a nonincreasing null sequence in $(0,1)$ and suppose that $(\beta_n)$ is a sequence of ordinals such that  $\sup \beta_m = \omega^\xi > \beta_n > 0$ for all $n$, $0 < \xi < \omega_1$.  Let 
\[ \gamma(\ep,m) = \max\{\beta_{n_s}+\cdots+\beta_{n_1}: \ep\theta_{n_s}\cdots\theta_{n_1} > \theta_m\}\ (\max \emptyset = 0).\]
The following are equivalent for any block subspace $Y$ of $T[(\theta_n,\cS_{\beta_n})^\infty_{n=1}]$.  
\begin{enumerate}
\item There exists $\ep > 0$ such that for all $\beta< \omega^\xi$, there exists $m \in \N$ such that $\gamma(\ep,m) + 2+\beta < \beta_m$ and every block subspace $Z$ of $Y$ contains a block sequence equivalent to a subsequence of $(e_k)$.
\item Every block subspace $Z$ of $Y$ contains an $\ell^1$-$\cS_{\omega^\xi}$-spreading model.
\item The Bourgain $\ell^1$-index $I_b(Z) = I(Z) > \omega^{\omega^\xi}$ for any block subspace $Z$ of $Y$.
\end{enumerate}
If one (and hence all) of these conditions holds for a block subspace $Y$ of $T[(\theta_n,\cS_{\beta_n})^\infty_{n=1}]$, then $Y$ is arbitrarily distortable.
Moreover, all these equivalent conditions hold for the space $T[(\theta_n,\cS_{\beta_n})^\infty_{n=1}]$
if $\xi $ is a limit ordinal.
\end{theorem}

When considering the mixed Tsirelson space $T[(\theta_n,\cS_{n})^\infty_{n=1}]$, it is customary to assume without loss of generality that $\theta_{m+n} \geq \theta_m\theta_n$ for all $m, n$.  In this case, it was shown in the proof of Corollary 28 in \cite{LT} that condition $(\dagger)$ is equivalent to $\lim_m\limsup_n\theta_{m+n}/\theta_n > 0$.

\begin{corollary}\label{corollary 12}
Let $(\theta_n)$ be a nonincreasing null sequence in $(0,1)$ such that $\theta_{m+n} \geq \theta_m\theta_n$ for all $m, n$.  
The following are equivalent for any block subspace $Y$ of $T[(\theta_n,\cS_{n})^\infty_{n=1}]$.  
\begin{enumerate}
\item $\lim_m\limsup_n\theta_{m+n}/\theta_n > 0$ and every block subspace $Z$ of $Y$ contains a block sequence equivalent to a subsequence of $(e_k)$.
\item Every block subspace $Z$ of $Y$ contains an $\ell^1$-$\cS_{\omega}$-spreading model.
\item The Bourgain $\ell^1$-index $I_b(Z) = I(Z) > \omega^{\omega}$ for any block subspace $Z$ of $Y$.
\end{enumerate}
If one (and hence all) of these conditions holds for a block subspace $Y$ of $T[(\theta_n,\cS_{n})^\infty_{n=1}]$, then $Y$ is arbitrarily distortable.
Moreover, all these equivalent conditions hold for the space $T[(\theta_n,\cS_{n})^\infty_{n=1}]$
if $\lim\theta_n^{1/n} = 1$.
\end{corollary}

\noindent{Remark.} The fact that every block subspace of $T[(\theta_n,\cS_{n})^\infty_{n=1}]$ contains an $\ell^1$-$\cS_{\omega}$-spreading model
if $\lim\theta_n^{1/n} = 1$ is due to Argyros, Deliyanni and Manoussakis \cite[Proposition 3.1]{ADM}.  Androulakis and Odell \cite{AO} showed that if $\lim \theta_n/\theta^n = 0$, where $\theta = \lim \theta_n^{1/n}$, then $T[(\theta_n,\cS_{n})^\infty_{n=1}]$ is arbitrarily distortable.

\begin{proof}
It suffices to prove the ``moreover" statement.  Clearly every normalized block sequence in $T[(\theta_n,\cS_{n})^\infty_{n=1}]$ is an $\ell^1$-$\cS_n$-spreading model with constant $\theta_n^{-1}$ for any $n$.  Thus for any block subspace $Y$, $I_b(Y,\theta_{m2^n}^{-1}) \geq \omega^{m2^n}$ for all $m, n$.
By the proof of Theorem 1.1 in \cite{JO}, it follows that $I_b(Y,\theta_{m2^n}^{-1/2^n}) \geq \omega^{m}$.  Using the hypothesis, we see that $I_b(Y,1+\ep) \geq \omega^{m}$ for all $\ep > 0$ and all $m$.  Hence $I_b(Y,1+\ep) \geq \omega^\omega$. By \cite[Lemma 5.7]{JO}, $I_b(Y) > \omega^\omega$.  This proves condition (3).
\end{proof}

For an ordinal $\beta$ with Cantor normal form $\beta = \omega^{\beta_1}\cdot m_1 + \dots + \omega^{\beta_n}\cdot m_n$, call $m_1$ the leading coefficient of $\beta$.  The preceding proof shows that if $(\F_n)$ is an increasing sequence of regular families so that $(\io(\F_n))$ increases nontrivially to $\omega^{\omega^{\xi}}$, where $\xi$ is a countable successor ordinal, and $\sup_n\theta_n^{{1}/{k_n}} = 1$, where $k_n$ is the leading coefficient of $\ell(\io(\F_n))$, then every block subspace of $T[(\theta_n,\F_{n})^\infty_{n=1}]$ contains an $\ell^1$-$\cS_{\omega^{\xi}}$-spreading model.

\begin{lemma} \label{lemma 4}
Let $\F$ be a regular family and let $M \in [\N]$.
\begin{enumerate}
\item If $0 < \io(\F) < \omega$, then there exists $N\in [M]$ such that $\F \cap [N]^{<\infty} = \A_j \cap [N]^{<\infty}$, where $j = \io(\F)$.
\item If $\io(\F) \geq \omega$, then there exists $N \in [M]$ such that $\F[\A_3] \cap [N]^{<\infty} \subseteq (\F)^2$.
\item If $\io(\F) \geq \omega$, then there exist $N \in [M]$ and $j \in \N$ such that $(\F\ominus\A_j)[\A_2] \cap [N]^{<\infty} \subseteq \F$.
\end{enumerate}
\end{lemma}

\begin{proof}
1. Since $\io(\F) = j$, $\F \subseteq \A_j$. On the other hand, choose $n_0 \in M$ such that $\{n_0\} \in \F^{(j-1)}$.  If $j = 1$, then set $N = \{n_0, n_0+1,\dots\}\cap M$.  Clearly $\A_1\cap [N]^{<\infty} \subseteq \F$.  If $j > 1$, consider $\G = \{G : n_0 < G, \{n_0\}\cup G \in \F\}$.  Then $\io(\G) = j-1$.  Using induction, we obtain $N_1 \in [M]$, $n_0 < N_1$ such that $\G\cap [N_1]^{<\infty} = \A_{j-1}\cap [N_1]^{<\infty}$.  Let $N = \{n_0\} \cup N_1$.  If $F \in \A_{j}\cap [N]^{<\infty}$, then $F$ is a spreading of $\{n_0\}\cup G$ for some $G \in \A_{j-1}\cap [N_1]^{<\infty}= \G\cap [N_1]^{<\infty}$.  Hence $F \in \F$.\\
2.  This follows from \cite[Theorem 1.1]{G} since $\io(\F[\A_3]) \leq \io(\A_3)\cdot\io(\F) < \io(\F)\cdot 2 = \io((\F)^2)$.\\
3.  Write $\io(\F) = \alpha + (j - 1)$ for some limit ordinal $\alpha$ and some $j \in \N$.  It is readily verified that $\F\ominus\A_j$ is a regular family and that $(\F\ominus\A_j)^{(\beta)} \subseteq \F^{(\beta)}\ominus\A_j$ for any $\beta$.  If $\io(\F\ominus\A_j) \geq \alpha$, then $\emptyset \in \F^{(\alpha)}\ominus\A_j$.  Hence there exists $A$, $|A| = j$ such that $A \in \F^{(\alpha)}$.  But then $\io(\F) \geq \alpha + j$, a contradiction.  Thus $\io(\F\ominus\A_j) < \alpha$.  It follows that $\io((\F\ominus\A_j)[\A_2]) < 2\cdot \alpha = \alpha \leq \io(\F)$.  By \cite[Theorem 1.1]{G}, there exists $N \in [M]$ such that $(\F\ominus\A_j)[\A_2] \cap [N]^{<\infty} \subseteq \F$.
\end{proof}

Given a regular family $\F$ and $M = (p_k) \in [\N]$, define the family $\MF$ by
$\MF = \{F: (p_k)_{k\in F} \in \F\}$.
It is clear that $\MF$ is a regular family.  Furthermore, the subspace $[(e_k)_{k\in M}]$ of 
$T[(\theta_n,\F_{n})^\infty_{n=1}]$ is easily seen to coincide with the mixed Tsirelson space 
$T[(\theta_n,\MF_{n})^\infty_{n=1}]$ under a natural identification.  The next proposition shows that the tameness of the sequence of regular families is not a restriction if one is allowed to pass to a subsequence of the unit vector basis.

\begin{proposition}
There exists $M \in [\N]$ and a tame sequence of regular families $(\G_n)$ such that $T[(\theta_n,\MF_{n})^\infty_{n=1}]$ is isomorphic to $T[(\theta_n,\G_{n})^\infty_{n=1}]$ via
the formal identity.
\end{proposition}

\begin{proof}
Let $m_0$ be the largest number such that $\alpha_{m_0} \leq \omega$.  (Take $m_0$ to be $0$ if $\alpha_n > \omega$ for all $n$.) Choose a strictly increasing sequence $(m_k)^\infty_{k=1}$ such that $m_1 > m_0$ and $\theta_{m_{k+1}} \leq \theta_{m_k}/2$ for all $k \in \N$.  By (1) and (2) of Lemma \ref{lemma 4}, there exists $M_0 \in [\N]$ such that for each $n \leq m_0$, either $\F_n \cap [M_0]^{<\infty} = \A_j\cap [M_0]^{<\infty}$ for some $j$, or $\F_n[\A_3] \cap [M_0]^{<\infty} \subseteq (\F_n)^2$. It is possible to choose a decreasing sequence $(M_k)^\infty_{k=1}$ of infinite subsets of $M_0$ and a sequence $(r_k)^\infty_{k=1}$ in $\N$ so that whenever $m_{k-1} < n \leq m_k$, $k \in \N$,
\begin{enumerate}
\item $\cS_1\cap[M_k]^{<\infty} \subseteq \F_n$ -- by \cite[Theorem 1.1]{G} since $\io(\cS_1) = \omega < \io(\F_n)$,
\item $\F_n[\A_3] \cap [M_k]^{<\infty} \subseteq (\F_n)^2$ -- by (2) of Lemma \ref{lemma 4},
\item $(\F_n\ominus\A_{r_k})[\A_2] \cap [M_k]^{<\infty} \subseteq \F_n$ -- by (3) of Lemma \ref{lemma 4}.
\end{enumerate}
Choose a strictly increasing sequence $(p_k)^\infty_{k=1}$ so that $r_k \leq p_k \in M_k$ 
for all $k \in \N$.
Define $M = (p_k)$ and set $\G_n = \MF_n$ if $n \leq m_0$ and $\G_n = \{G \in \MF_n: G \geq k\}$ if $m_{k-1} < n \leq m_{k}$, $k \in \N$.  By \cite[Proposition 1]{LT}, 
$T[(\theta_n,\MF_{n})^\infty_{n=1}]$ is isomorphic to $T[(\theta_n,\G_{n})^\infty_{n=1}]$ via
the formal identity.  It remains to show that the sequence $(\G_n)$ is tame.

First suppose that $n \leq m_0$.  If $\F_n \cap [M_0]^{<\infty} = \A_j\cap[M_0]^{<\infty}$ for some $j$, then clearly $\G_n = \MF_n = \A_j$.  Otherwise, $\F_n[\A_3]\cap [M]^{<\infty} \subseteq \F_n[\A_3]\cap [M_0]^{<\infty} \subseteq (\F_n)^2$.  If $G \in \G_n[\A_3] = \MF_n[\A_3]$, then $(p_k)_{k\in G} \in \F_n[\A_3]\cap [M]^{<\infty} \subseteq (\F_n)^2$.  Hence $G \in (\MF_n)^2$.

Now assume that $n > m_0$.  Choose $k$ such that $m_{k-1} < n \leq m_k$.  If $G \in \G_n[\A_3]$, then $G \in \MF_n[\A_3]$ and $G \geq k$.  Hence $p_k \leq (p_i)_{i\in G} \in \F_n[\A_3]$.  Thus $(p_i)_{i\in G} \in \F_n[\A_3]\cap[M_k]^{<\infty}\subseteq (\F_n)^2$.  Therefore $G \in (\MF_n)^2$ and $G \geq k$.  It follows that $G \in (\G_n)^2$.

Finally, we show that $(\G_n\ominus\G_m)[\A_2] \subseteq \G_n$ whenever $n > m > m_0$.  
Choose $k$ and $l$ such that  and $m_{k-1} < n \leq m_k$ and $m_{l-1} < m \leq m_l$.
Suppose that $G \in \G_n\ominus\G_m$.  There is a maximal $H \in \G_m$ such that $H < G$ and $H \cup G \in \G_n$.  
We claim that $|H| \geq r_k$.  Indeed,
by definition of $\G_n$, $H \geq k$.  Thus $r_k \leq p_k \leq (p_i)_{i\in H}$.  If $|H| < r_k$, there exists a nonempty set $I > H$ such that $(p_i)_{i\in H\cup I} \in \cS_1$.  Clearly, $(p_i)_{i\in H\cup I} \in [M_k]^{<\infty} \subseteq [M_l]^{<\infty}$ as well.  Therefore, $(p_i)_{i\in H\cup I} \in \F_m$ by condition (1) above.
By definition, $H \cup I \in \MF_m$.  Since $H\cup I \geq k \geq l$, $H\cup I \in \G_m$, contrary to the maximality of $H$.  This proves the claim.
It follows from the claim that $(p_i)_{i\in G} \in \F_n\ominus\A_{r_k}$.  Thus $(p_i)_{i\in J} \in (\F_n\ominus\A_{r_k})[\A_2]$ for all $J \in (\G_n\ominus\G_m)[\A_2]$.
Clearly, for such $J$, $J \geq k$ and hence $(p_i)_{i\in J}$ is  in $[M_k]^{<\infty}$.  Therefore,
\[ (p_i)_{i\in J} \in (\F_n\ominus\A_{r_k})[\A_2] \cap [M_k]^{<\infty} \subseteq \F_n \]
by condition (3) above.  This shows that $J \in \MF_n$.  As $J \geq k$, we have $J \in \G_n$, as desired.
\end{proof}

\begin{corollary}
Suppose that either {\em (a)} $\xi$ is a countable limit ordinal or that {\em (b)} $\xi$ is a countable successor ordinal and $\sup_n\theta_n^{1/k_n} = 1$, where $k_n$ is the leading coefficient of $\ell(\alpha_n)$.  Then there exists $M \in [\N]$ such that the subspace $Y = [(e_k)_{k\in M}]$ of $X$ has the following properties.
\begin{enumerate}
\item Every block subspace of $Y$ has an $\ell^1$-$\cS_{\omega^\xi}$-spreading model.
\item Every block subspace of $Y$ contains a block sequence equivalent to a subsequence of $(e_k)_{k\in M}$.
\item $Y$ is arbitrarily distortable.
\end{enumerate}
\end{corollary}

Schlumprecht proposed a classification of Banach spaces as follows \cite{S2}.  A Banach space with a normalized basis $(u_k)$ is said to be {\em Class $1$} if every normalized block sequence has a subsequence equivalent to a subsequence of $(u_k)$.
It is {\em Class $2$} if every block subspace contains two block sequences $(y_k)$ and $(z_k)$ so that the map $y_k \mapsto z_k$ extends to a bounded linear strictly singular operator.  Recall that an operator is {\em strictly singular} if its restriction to any infinite dimensional subspace is not an isomorphism.  Schlumprecht asks whether every infinite dimensional Banach space contains a subspace with a basis that is either Class $1$ or Class $2$.  He also proved a criterion for a Banach space to be Class $2$ \cite[Theorem 1.4 and Corollary 1.5]{S2}.  We conclude with a note showing that his proof applies to mixed Tsirelson spaces satisfying the hypotheses of Theorem \ref{theorem 8}.  A Banach space is $c_0$-{\em saturated} if every closed infinite dimensional subspace contains an isomorphic copy of $c_0$.

\begin{proposition}
Let $Y$ be a block subspace of a mixed Tsirelson space $X$ and suppose that $Y$ satisfies all the hypotheses of Theorem \ref{theorem 8}.  Then $Y$ is a Class $2$ space.
\end{proposition}

\begin{proof}
Denote by $(e_k)$ the unit vector basis of $X$. We will show below that there are a regular family $\G$ with $\io(\G) \leq  \omega^{\omega^\xi}$ and a finite constant $C$ so that $\|\sum a_ke_k\| \leq C\sup_{G\in \G}\sum_{k\in G}|a_k|$ for all $(a_k) \in c_{00}$.  Denote the unit vector basis in $c_{00}$ by $(u_k)$ and let $U$ be the completion of $c_{00}$ with respect to the norm $\|\sum a_ku_k\| = \sup_{G\in \G}\sum_{k\in G}|a_k|$ for all $(a_k) \in c_{00}$.  The map that sends $\sum a_ku_k$ to the function on $\G$ given by $G \mapsto \sum_{k\in G}a_k$ is an embedding of $U$ into $C(\G)$, the space of continuous functions on the countable compact metric space $\G$.  Hence $U$ is $c_0$-saturated.
Let $Z$ be a block subspace of $Y$.  By the hypothesis, there is a block sequence $(z_k)$ in $Z$ that is equivalent to a subsequence $(e_{m_k})$ of $(e_k)$.  Also, there is a sequence $(y_k)$ in $Z$ that generates an $\ell^1$-${\cS_{\omega^\xi}}$-spreading model.  We may replace $(y_k)$ with an appropriate subsequence of $(y_{2k}-y_{2k+1})$ if necessary to assume that $(y_k)$ is equivalent to a block sequence.  By definition of the norm in $X$, there is a positive constant $K$ so that $\|\sum a_ky_k\| \geq K^{-1}\sum_{k\in F}|a_k|$ for all $F \in \F_1[\cS_{\omega^{\xi}}]$.
Since $\io(\F_1) > 1$ by assumption, $\io(\F_1[\cS_{\omega^{\xi}}]) > \omega^{\omega^{\xi}} \geq \io(\G)$.
Using \cite[Theorem 1.1]{G} and replacing $M = (m_k)$ with a subsequence if necessary, we may assume that $\G\cap [M]^{<\infty} \subseteq \F_1[\cS_{\omega^\xi}]$.  Because $(z_k)$ is equivalent to $(e_{m_k})$ and $(y_k)$ is equivalent to a block sequence, it follows that the map $y_{m_k} \mapsto z_k$ extends to a bounded linear map $T: [(y_{m_k})]\to [(z_k)]$. 
Now, for all $(a_k) \in c_{00}$,
\begin{align*}
\|\sum a_ku_{m_k}\| &= \sup_{G\in \G}\sum_{m_k\in G}|a_k| \\
&\leq \sup_{G\in\F_1[\cS_{\omega^{\xi}}]}\sum_{m_k\in G}|a_k| \\
&\leq K\|\sum a_ky_{m_k}\|.
\end{align*}
Hence $y_{m_k}\mapsto u_{m_k}$ extends to a bounded linear map $S: [(y_{m_k})]\to [(u_{m_k})]$.
However, $(z_k)$ is equivalent to $(e_{m_k})$ and 
\[ \|\sum a_ke_{m_k}\| \leq C\sup_{G\in \G}\sum_{m_k\in G}|a_k| = C\|\sum a_ku_{m_k}\|. \]
Thus $u_{m_k} \mapsto z_k$ extends to a bounded linear map $R: [(u_{m_k})] \to [(z_k)]$.  Therefore, $T = RS$ is a factorization of $T$ through the $c_0$-saturated space $[(u_{m_k})]$.  Since $[(y_{m_k})]$ does not contain a copy of $c_0$, $T$ is strictly singular.

It remains to show the existence of the family $\G$.  Choose a strictly increasing sequence $(n_i)$ such that $\pi_{i} < 2^{-i}$ for all $i$, where 
\[ \pi_i = \max\{\theta_{m_1}\cdots\theta_{m_r}: m_1 + \dots + m_r > n_i\}. \]
For each $i$, let $\G_i = \cup\{[\F_{m_r},\dots,\F_{m_1}]: m_1+\dots+m_r \leq n_i\}$.  Here $[\F_{m_r},\dots,\F_{m_1}]$ is defined inductively as $\F_{m_r}[\F_{m_{r-1}},\dots,\F_{m_1}]$.  It follows from \cite[Proposition 12]{LT1} that $\io(\G_i) < \omega^{\omega^\xi}$ since $\io(\F_n) < \omega^{\omega^\xi}$ for each $n$.
Let $\G$ consist of all sets $G$ such that $G \in \G_i$ for some $i \leq G$ together with all singletons.  Then $\io(\G) \leq \omega^{\omega^\xi}$.  For any $x = \sum a_ke_k$, $(a_k) \in c_{00}$, let $\T$ be an admissible tree that norms $x$.  Denote by $\E$ the set of all leaves of $\T$.  Also, if $t(E) = \theta_{m_1}\cdots\theta_{m_r}$, $E \in \E$, set $r(E) = m_1+\dots+m_r$.  Note that $\{E\in \E: r(E)\leq n_i\}$ is $\G_i$-admissible.  Thus
\begin{align*}
\T x &= \sum_{E\in \E}t(E)\|Ex\|_{c_0}\\
& = \sum^\infty_{i=1}\sum_{n_{i-1}<r(E)\leq n_i}t(E)\|Ex\|_{c_0}\\
&\leq \sum^\infty_{i=1}\pi_{{i-1}}\rho_i(x),
\end{align*}
where $\rho_i(x) = \sup_{G\in \G_i}\sum_{k\in G}|a_k|$.  However,
\[ \rho_i(x) \leq \sum^i_{k=1}|a_k| + \sup_{G\in \G_i}\sum_{k\in G, k > i}|a_k| \leq i\|x\|_{c_0}+ \sup_{G\in \G}\sum_{k\in G}|a_k|.\]
Therefore, 
\begin{align*}
\|x\| &\leq \sum^\infty_{i=1}\pi_{{i-1}}\rho_i(x) \leq \sum^\infty_{i=1}\frac{\rho_i(x)}{2^{i-1}} \\
&\leq \|x\|_{c_0}\sum^\infty_{i=1}\frac{i}{2^{i-1}} +
\sum^\infty_{i=1}\frac{1}{2^{i-1}}\sup_{G\in \G}\sum_{k\in G}|a_k|\\
&\leq 6\sup_{G\in \G}\sum_{k\in G}|a_k|.
\end{align*}
This completes the proof.
\end{proof}

\end{document}